\documentclass[10pt]{amsart}
\usepackage{amscd}
\usepackage{amssymb}
\usepackage[all]{xy}
\usepackage{graphicx}
\author{Amnon Yekutieli}
\title[An Averaging Process]
{An Averaging Process for Unipotent Group Actions} 
\date{12 January 2006}
\address{Department of  Mathematics,
Ben Gurion University, Be'er Sheva 84105, Israel}
\email{amyekut@math.bgu.ac.il}
\thanks{
Primary: 14L30; Secondary: 18G30, 20G15.}
\thanks{ This work was partially supported by the 
US - Israel Binational Science Foundation}

\newtheorem{thm}[equation]{Theorem}
\newtheorem{cor}[equation]{Corollary}

\newtheorem{lem}[equation]{Lemma}
\theoremstyle{definition}
\newtheorem{dfn}[equation]{Definition}
\newtheorem{rem}[equation]{Remark}

\newtheorem{sce}[equation]{Scenario}

\numberwithin{equation}{section}

\newcommand{\opn}{\operatorname}
\newcommand{\cat}[1]{\operatorname{\mathsf{#1}}}

\newcommand{\ul}{\underline}

\newcommand{\rmitem}[1]{\item[\text{\textup{(#1)}}]}
\newcommand{\mfrak}[1]{\mathfrak{#1}}

\newcommand{\mbf}[1]{\mathbf{#1}}
\newcommand{\mrm}[1]{\mathrm{#1}}
\newcommand{\mbb}[1]{\mathbb{#1}}

\newcommand{\smfrac}[2]{{\textstyle \frac{#1}{#2}}}
\newcommand{\tup}[1]{\textup{#1}}
\newcommand{\bsym}[1]{\boldsymbol{#1}}

\newcommand{\til}[1]{\tilde{#1}}

\newcommand{\set}[1]{\{ #1 \}}

\newcommand{\K}{\mbb{K}}

\begin{document}

\begin{abstract}
We present an averaging process for sections of a torsor under a 
unipotent group. This process allows one to integrate local sections 
of such a torsor into a global simplicial section. 
The results of this paper have 
applications to deformation quantization of algebraic varieties. 
\end{abstract}

\maketitle

\setcounter{section}{-1}
\section{Introduction}

Let $\K$ be a field of characteristic $0$. For any natural number
$q$ let $\K [t_0, \ldots, t_q]$ be the polynomial algebra, and let
$\bsym{\Delta}^q_{\K}$ be the {\em geometric 
$q$-dimensional simplex}
\[ \bsym{\Delta}^q_{\K} := \opn{Spec}
\K [t_0, \ldots, t_q] / (t_0 + \cdots + t_q - 1) . \] 
Our main result is the following theorem. 

\begin{thm} \label{thm0.1}
Let $G$ be a unipotent algebraic group over 
$\K$, let $X$ be a $\K$-scheme, let $Z$ be a $G$-torsor over $X$, 
and let $Y$ be any $X$-scheme. Suppose
$\bsym{f} = (f_0, \ldots, f_q)$ is a sequence of $X$-morphisms
$f_i : Y \to Z$. Then there is an $X$-morphism
\[ \opn{wav}_G(\bsym{f}) : \bsym{\Delta}^q_{\K} \times Y \to Z  \]
called the {\em weighted average}, such that 
the operation $\opn{wav}_G$ is 
symmetric, simplicial, functorial in the data $(G, X, Y, Z)$, 
and is the identity for $q = 0$. 
\end{thm}

``Symmetric'' means that $\opn{wav}_G$ is equivariant for the 
action of the permutation group of $\{0, \ldots, q \}$ on the 
sequence $\bsym{f}$ and the scheme $\bsym{\Delta}^q_{\K}$.
``Simplicial'' says that as $q$ varies
\[ \opn{wav}_G : \opn{Hom}_X(Y, Z)^{q+1} \to
\opn{Hom}_X(\bsym{\Delta}^q_{\K} \times Y, Z) \]
is a map of simplicial sets. ``Functorial in the data 
$(G, X, Y, Z)$'' means that e.g.\ given a morphism 
$h : Z \to Z'$ of $G$-torsors over $X$, one has
\[ h \circ \opn{wav}_G(\bsym{f}) = \opn{wav}_G(h \circ \bsym{f}), 
\]
where $h \circ \bsym{f}$ is the sequence
$(h \circ f_0, \ldots, h \circ f_q)$. 
Theorem \ref{thm0.1} is repeated, in full detail, as 
Theorem \ref{thm5.1}. 

Observe that when we restrict the morphism 
$\opn{wav}_G(\bsym{f})$ to each of the vertices of 
$\bsym{\Delta}^q_{\K}$ we recover the original morphisms 
$f_0, \ldots, f_q$. This is due to the simplicial property of 
$\opn{wav}_G$. Thus $\opn{wav}_G(\bsym{f})$ interpolates 
between $f_0, \ldots, f_q$. This is illustrated in Figure 
\ref{fig1}.

Here is a naive corollary of Theorem \ref{thm0.1},
which can further explain the result. Let us 
write $G(\K)$ for the group of $\K$-rational points of $G$. 
By a {\em weight sequence} in $\K$ we mean a sequence 
$\bsym{w} = (w_0, \ldots, w_q)$ of elements of $\K$ such that
$\sum w_i = 1$. 

\begin{cor} \label{cor0.1}
Let $G$ be a unipotent group over $\K$.
Suppose $Z$ is a set with $G(\K)$-action which is transitive and 
has trivial stabilizers. Let 
$\bsym{z} = (z_0, \ldots, z_q)$ be a sequence of points in $Z$ and 
let $\bsym{w}$ be a weight sequence in $\K$. Then there is a 
point 
$\opn{wav}_{G, \bsym{w}}(\bsym{z}) \in Z$
called the weighted average. 
The operation $\opn{wav}_{G}$
is symmetric, functorial, simplicial, 
and is the identity for $q = 0$. 
\end{cor}

The corollary is proved in Section 1. The idea is of course to 
take $Y = X := \opn{Spec} \K$ in the theorem, and to note that 
$\bsym{w}$ is a $\K$-rational point of $\bsym{\Delta}^q_{\K}$. 

Another consequence of Theorem \ref{thm0.1} is a new proof of the 
fact that a unipotent group in characteristic $0$ is special. See 
Remark \ref{rem1.1}. (This observation is due to Reichstein.)

Theorem \ref{thm0.1} was discovered in the course of research on 
deformation quantization in algebraic geometry \cite{Ye}, in 
which we tried to apply ideas of Kontsevich \cite{Ko}
to the algebraic context. Here is a brief outline. Suppose $X$ is 
a smooth $n$-dimensional $\K$-scheme. The coordinate bundle of $X$ 
\cite{GK, Ko} is an infinite dimensional bundle $Z \to X$ which 
parameterizes formal coordinate systems on $X$. The bundle $Z$ is 
a torsor under an affine group scheme 
$\mrm{GL}_{n, \K} \ltimes G$, where $G$ is pro-unipotent. One is 
interested in sections of the quotient bundle 
$\bar{Z} := Z / \mrm{GL}_{n, \K}$. If we are in the 
differentiable setup (i.e.\ $\K = \mbb{R}$ and $X$ is a 
$\mrm{C}^{\infty}$ manifold) then the fibers of $\bar{Z}$ are 
contractible (since they are isomorphic to $G(\mbb{R})$), 
and therefore global $\mrm{C}^{\infty}$ sections exist. 
However in the setup of algebraic geometry such an argument is 
invalid, and we were forced to seek an alternative approach.

Our solution was to use {\em simplicial sections} of $\bar{Z}$
(see Section 2, and especially Corollary \ref{cor2.1}). This 
approach was inspired by work of Bott \cite{Bo} on simplicial 
connections (cf.\ also \cite{HY}).

\medskip \noindent
\textbf{Acknowledgments.} 
The author thanks David Kazhdan and Zinovy Reichstein for useful 
conversations. Also thanks to the referee for reading the paper 
carefully and suggesting a few improvements.

\section{The Averaging Process}
\label{sec1}

Throughout the paper $\K$ denotes a fixed base field of 
characteristic $0$. All schemes and all morphisms are over $\K$.

We begin by recalling some standard facts about the combinatorics
of simplicial objects. 
Let $\bsym{\Delta}$ denote the category with objects the ordered 
sets $[q] := \{ 0, 1, \ldots, q \}$, $q \in \mbb{N}$. The 
morphisms $[p] \to [q]$ are the order preserving functions,
and we write
$\bsym{\Delta}^q_p := \opn{Hom}_{\bsym{\Delta}}([p], [q])$. 
The $i$-th co-face map 
$\partial^i : [p] \to [p + 1]$
is the injective function that does not take the value $i$;
and the $i$-th co-degeneracy map
$\mrm{s}^i : [p] \to [p - 1]$
is the surjective function that takes the value $i$ twice.
All morphisms in $\bsym{\Delta}$ are compositions of various 
$\partial^i$ and $\mrm{s}^i$. 

An element of $\bsym{\Delta}^q_p$ may be thought of as a sequence
$\bsym{i} = (i_0, \ldots, i_p)$ of integers with
$0 \leq i_0 \leq \cdots \leq i_p \leq q$. 
Given $\bsym{i} \in \bsym{\Delta}^m_q$, 
$\bsym{j} \in \bsym{\Delta}_m^p$ and
$\alpha \in \bsym{\Delta}^q_p$, 
we sometimes write 
$\alpha_*(\bsym{i}) := \bsym{i} \circ \alpha \in \bsym{\Delta}^m_p$ 
and 
$\alpha^*(\bsym{j}) := \alpha  \circ\bsym{j} \in 
\bsym{\Delta}_m^q$.

Let $\cat{C}$ be some category. A {\em cosimplicial object} in 
$\cat{C}$ is a functor $C : \bsym{\Delta} \to \cat{C}$.
We shall usually refer to the cosimplicial object as
$C = \set{C^p}_{p \in \mbb{N}}$,
and for any $\alpha \in \bsym{\Delta}_p^q$ the corresponding 
morphism in $\cat{C}$ will be denoted by 
$\alpha^* : C^p \to C^q$.
A {\em simplicial object} in 
$\cat{C}$ is a functor $C : \bsym{\Delta}^{\mrm{op}} \to \cat{C}$.
The notation for a simplicial object
will be $C = \set{C_p}_{p \in \mbb{N}}$ and
$\alpha_* : C_q \to C_p$.

An important example is the cosimplicial scheme 
$\{ \bsym{\Delta}^p_{\K} \}_{p \in \mbb{N}}$.
The morphisms are defined as follows. 
For any $p$ we identify the ordered set 
$[p]$ with the set of vertices of $\bsym{\Delta}^p_{\K}$. 
Given $\alpha \in \bsym{\Delta}^q_{p}$ the morphism
$\alpha^* : \bsym{\Delta}^p_{\K} \to \bsym{\Delta}^q_{\K}$
is then the unique linear morphism extending 
$\alpha : [p] \to [q]$. 

Let $G$ be a unipotent (affine finite type)
algebraic group over $\K$, with 
(nilpotent) Lie algebra $\mfrak{g}$. We write $d(G)$ for the minimal 
number $d$ such that there exists a chain of closed normal subgroups
$G = G_0 \supset G_1 \cdots \supset G_d = 1$
with $G_k / G_{k+1}$ abelian 
for all $k \in \{ 0, \ldots, d-1 \}$. The exponential map
$\opn{exp}_G : \mfrak{g} \to G$ is an isomorphism
of schemes, with inverse $\opn{log}_G$; see 
\cite[Theorem VIII.1.1]{Ho}.

Given a $\K$-scheme $X$ there is a Lie algebra 
$\mfrak{g} \times X$ (in the category of $X$-schemes), and a 
group-scheme $G \times X$. There is also an induced exponential 
map
\[ \opn{exp}_{G \times X} :=
\opn{exp}_G \times \bsym{1}_X  : 
\mfrak{g} \times X \to G \times X . \] 
We will need the following result. 

\begin{lem} \label{lem1.1}
Let $G, G'$ be two unipotent groups, with Lie algebras
$\mfrak{g}, \mfrak{g}'$. Let $X, X'$ be schemes, let 
$X \to X'$ be a morphism of schemes, and let
$\phi : G \times X \to G' \times X'$ be a morphism of 
group-schemes over $X'$. Denote by 
$\mrm{d} \phi : \mfrak{g} \times X \to \mfrak{g}' \times X'$
the induced Lie algebra morphism \tup{(}the differential of 
$\phi$\tup{)}. Then the diagram
\begin{equation} \label{eqn1.6}
\begin{CD}
\mfrak{g} \times X @>{\opn{exp}_{G \times X}}>> G \times X \\
@V{\mrm{d} \phi}VV @V{\phi}VV \\
\mfrak{g}' \times X' @>{\opn{exp}_{G' \times X'}}>> G' \times X'
\end{CD} 
\end{equation}
commutes.
\end{lem}

\begin{proof}
For the case $X = X' = \opn{Spec} \K$ this is contained in the 
proof of \cite[Theorem VIII.1.2]{Ho}.

In order to handle the general case we first recall the 
Campbell-Baker-Hausdorff formula:
\[ \opn{exp}_{G}(\gamma_1) \cdot \opn{exp}_{G}(\gamma_2) =
\opn{exp}_{G}(F(\gamma_1, \gamma_2)) , \]
where 
\[ F(\gamma_1, \gamma_2) = \gamma_1 + \gamma_2 + 
\smfrac{1}{2} [\gamma_1, \gamma_2] + \cdots \]
is a universal power series (see \cite[Section XVI.2]{Ho}).
Hence if we define 
$\gamma_1 * \gamma_2 := F(\gamma_1, \gamma_2)$,
then $(\mfrak{g}, *)$ becomes an algebraic group (with unit 
element $0$), and 
$\opn{exp}_{G} : (\mfrak{g}, *) \to (G, \cdot)$
is a group isomorphism. In this way we may eliminate $G$ 
altogether, and just look at the scheme $\mfrak{g}$ with its two 
structures: a Lie algebra and a group-scheme. Note that now
$\mfrak{g}$ is its own Lie algebra, as can be seen from the $2$-nd 
order term in the series $F(\gamma_1, \gamma_2)$. 

Consider a morphism 
$\phi : \mfrak{g} \times X \to \mfrak{g}' \times X'$
of $X'$-schemes. Then $\phi$ is a morphism of Lie algebras over 
$X'$ iff it is a morphism of group-schemes (for the 
multiplications $*$). And moreover $\mrm{d} \phi = \phi$. 
Therefore the diagram (\ref{eqn1.6}) is commutative.
\end{proof}

From now on we shall write $\opn{exp}_{G}$ instead of
$\opn{exp}_{G \times X}$, for the sake of brevity.

Consider the following setup: $X$ is a $\K$-scheme, 
and $Y, Z$ are two $X$-schemes.
Suppose $Z$ is a torsor under the group scheme $G \times X$. 
We denote the action of $G$ 
on $Z$ by $(g, z) \mapsto g \cdot z$.

Let $f_0, \ldots, f_q : \bsym{\Delta}^q_{\K} \times Y \to Z$ 
be $X$-morphisms. 
We are going to define a new sequence of $X$-morphisms 
$f'_0, \ldots, f'_q : \bsym{\Delta}^q_{\K} \times Y \to Z$.
Because $Z$ is a torsor under $G \times X$, for any 
$i, j \in \set{0, \ldots, q}$ there exists a unique morphism
$g_{i, j} : \bsym{\Delta}^q_{\K} \times Y \to G$ 
such that $f_j(w, y) = g_{i, j}(w, y) \cdot f_i(w, y)$
for all points $w \in \bsym{\Delta}^q_{\K}$ and 
$y \in Y$. Here $w$ and $y$ are scheme-theoretic points,
i.e.\ 
$w \in \bsym{\Delta}^q_{\K}(U) = 
\opn{Hom}_{\K}(U, \bsym{\Delta}^q_{\K})$ 
and 
$y \in Y(U) = \opn{Hom}_{\K}(U, Y)$ 
for some $\K$-scheme $U$. Define
\begin{equation} \label{eqn6.X1}
f'_i(w, y) := \opn{exp}_G \Big( \sum_{j=0}^q t_j(w) \cdot
\opn{log}_G(g_{i, j}(w, y)) \Big) \cdot f_i(w, y) . 
\end{equation}
In this formula we view $t_j$ as a morphism 
$t_j : \bsym{\Delta}^q_{\K} \to \mbf{A}^1_{\K}$,
and we use the fact that $\mfrak{g}$ is a vector space
(in the category of $\K$-schemes).

For any set $S$ let us write 
$S^{\bsym{\Delta}^q_{0}}$ for the set of functions 
$\bsym{\Delta}^q_{0} \to S$, which is the same as the set of
sequences $(s_0, \ldots, s_q)$ in $S$. As usual 
$\opn{Hom}_X(\bsym{\Delta}^q_{\K} \times Y, Z)$
is the set of $X$-morphisms $\bsym{\Delta}^q_{\K} \times Y \to Z$. 
In this notation the sequences
$(f_0, \ldots, f_q)$ and $(f'_0, \ldots, f'_q)$ are
elements of 
$\opn{Hom}_X(\bsym{\Delta}^q_{\K} \times Y, Z)
^{\bsym{\Delta}^q_{0}}$. We denote by
\[ \opn{wsym}_G : \opn{Hom}_X(\bsym{\Delta}^q_{\K} \times Y, Z)
^{\bsym{\Delta}^q_{0}} \to
\opn{Hom}_X(\bsym{\Delta}^q_{\K} \times Y, Z)
^{\bsym{\Delta}^q_{0}} \]
the operation
$(f_0, \ldots, f_q) \mapsto (f'_0, \ldots, f'_q)$
given by the formula (\ref{eqn6.X1}). 

We will also need a similar operation
\[ \opn{w}_G : \opn{Hom}_X(Y, Z)^{\bsym{\Delta}^q_{0}} \to
\opn{Hom}_X(\bsym{\Delta}^q_{\K} \times Y, Z)
^{\bsym{\Delta}^q_{0}} , \]
defined by
\[ \opn{w}_G(f_0, \ldots, f_q) := (f'_0, \ldots, f'_q) \]
with
\begin{equation} \label{eqn1.3}
f'_i(w, y) := \opn{exp}_G \Big( \sum_{j=0}^q t_j(w) \cdot
\opn{log}_G(g_{i, j}(y)) \Big) \cdot f_i(y) .
\end{equation}

It is clear that for $q = 0$ both operations 
$\opn{w}_G$ and $\opn{wsym}_G$ act as the identity, i.e.\ 
$\opn{w}_G(f_0) = \opn{wsym}_G(f_0) = f_0$
for all $f_0 \in \opn{Hom}_X(Y, Z)$.
Both operations $\opn{w}_G$ and $\opn{wsym}_G$ are symmetric, 
namely they are equivariant for the simultaneous action of the 
permutation group of $\{ 0, \ldots, q \}$ on 
$\bsym{\Delta}^q_{0}$ and $\bsym{\Delta}^q_{\K}$.
Also if 
$\bsym{f} = (f_0, \ldots, f_q) \in 
\opn{Hom}_X(\bsym{\Delta}^q_{\K} \times Y, Z)
^{\bsym{\Delta}^q_{0}}$
is a constant sequence, i.e.\ $f_0 = \cdots = f_q$, then 
$\opn{wsym}_G(\bsym{f}) = \bsym{f}$.

\begin{lem} \label{lem6.X1}
Both operations $\opn{w}_G$ and $\opn{wsym}_G$ are simplicial.
Namely, given $\alpha \in \bsym{\Delta}^q_{p}$ the diagram
\[ \begin{CD}
\opn{Hom}_X(\bsym{\Delta}^q_{\K} \times Y, Z)
^{\bsym{\Delta}^q_{0}}
@>{\opn{wsym}_G}>>
\opn{Hom}_X(\bsym{\Delta}^q_{\K} \times Y, Z)
^{\bsym{\Delta}^q_{0}} \\
@V{\alpha_*}VV @V{\alpha_*}VV \\
\opn{Hom}_X(\bsym{\Delta}^p_{\K} \times Y, Z)
^{\bsym{\Delta}^p_{0}}
@>{\opn{wsym}_G}>>
\opn{Hom}_X(\bsym{\Delta}^p_{\K} \times Y, Z)
^{\bsym{\Delta}^p_{0}}
\end{CD} \]
is commutative, and likewise for $\opn{w}_G$.
\end{lem}

\begin{proof}
It suffices to consider $\alpha = \partial^i$ or 
$\alpha = \mrm{s}^i$. Since $\opn{wsym}_G$ is symmetric,
we may assume that 
$\alpha = \partial^q : [q-1] \to [q]$ or
$\alpha = \mrm{s}^q : [q+1] \to [q]$.
Fix a sequence
$\bsym{f} = (f_0, \ldots, f_q) \in 
\opn{Hom}_X(\bsym{\Delta}^q_{\K} \times Y, Z)
^{\bsym{\Delta}^q_{0}}$.
Let $g_{i,j} \in \opn{Hom}_{\K}(\bsym{\Delta}^q_{\K} \times Y, G)$
be such that
$f_{j} = g_{i,j} \cdot f_i$,
and let
$\bsym{f}' = (f'_0, \ldots, f'_q) :=  
\opn{wsym}_G(\bsym{f})$. 

First let's look at the case $\alpha = \partial^q$. 
Take $w \in \bsym{\Delta}^{q-1}_{\K}$ and
$y \in Y$, and let $v := \alpha^*(w) \in \bsym{\Delta}^q_{\K}$.
The coordinates of $v$ are
$t_j(v) = t_j(w)$ for $j \leq q-1$, and $t_q(v) = 0$.
Then for any $i$ the $i$-th term in the sequence
$\alpha_*(\bsym{f}')$, evaluated at $(w,y)$, equals
\begin{equation} \label{eqn6.X2}
\begin{aligned}
f'_i(v,y) & =
\opn{exp}_G \Big( \sum_{j=0}^q t_j(v) \cdot
\opn{log}_G(g_{i, j}(v, y)) \Big) \cdot f_i(v, y) \\
& = \opn{exp}_G \Big( \sum_{j=0}^{q-1} t_j(w) \cdot
\opn{log}_G(g_{i, j}(v, y)) \Big) \cdot f_i(v, y) .
\end{aligned}
\end{equation}
On the other hand, the $i$-th term of the 
sequence $\opn{wsym}_G(\alpha_*(\bsym{f}))$ is
\[ \begin{aligned}
& \opn{exp}_G \Big( \sum_{j=0}^{q-1} t_j(w) \cdot
\opn{log}_G( \alpha_*(g_{i, j})(w, y)) \Big) \cdot 
\alpha_*(f_i)(w, y) \\
& \qquad = \opn{exp}_G \Big( \sum_{j=0}^{q-1} t_j(w) \cdot
\opn{log}_G( g_{i, j}(v, y)) \Big) \cdot f_i(v, y) .
\end{aligned} \]
So indeed
$\alpha_* \circ \opn{wsym}_G = \opn{wsym}_G \circ\, \alpha_*$ in 
this case.

Next consider the case $\alpha = \mrm{s}^q$. Take
$w \in \bsym{\Delta}^{q+1}_{\K}$ and
$y \in Y$, and let $v := \alpha^*(w) \in \bsym{\Delta}^q_{\K}$.
The coordinates of $v$ are
$t_j(v) = t_j(w)$ for $j \leq q-1$, and 
$t_q(v) = t_q(w) + t_{q+1}(w)$.
For any $i \leq q$ the $i$-th term in the sequence 
$\alpha_*(\bsym{f})$, evaluated at $(w,y)$, is
$f'_i(v,y)$, which was calculated in (\ref{eqn6.X2}). 
The $(q+1)$-st term is also $f'_q(v,y)$.
On the other hand, 
for any $i \leq q$ the $i$-th term in the sequence 
$\opn{wsym}_G(\alpha_*(\bsym{f}'))$, evaluated at $(w,y)$, is
\[ z_i := \opn{exp}_G \Big( \sum_{j=0}^{q+1} t_j(w) \cdot
\opn{log}_G( \alpha_*(g_{i, j})(w, y)) \Big) \cdot 
\alpha_*(f_i)(w, y) . \]
But 
$t_q(w) + t_{q+1}(w) = t_q(v)$, 
$\alpha_*(g_{i,j})(w,y) = g_{i,j}(v,y)$ for $j \leq q$, and
\linebreak $\alpha_*(g_{i, q+1})(w, y) = g_{i, q}(v, y)$.
Therefore
\[ z_i = \opn{exp}_G \Big( \sum_{j=0}^{q} t_j(v) \cdot
\opn{log}_G( g_{i, j}(v, y)) \Big) \cdot f_i(v, y) . \]
For $i = q+1$ one has $z_{q+1} = z_q$. We conclude that 
$\alpha_* \circ \opn{wsym}_G = \opn{wsym}_G \circ\, \alpha_*$ in 
this case too.

The proof for $\opn{w}_G$ is the same.
\end{proof}

\begin{lem} \label{lem6.X2}
Both operations $\opn{w}_G$ and $\opn{wsym}_G$ are functorial in 
the data \linebreak $(G, X, Y, Z)$. Namely, given another such 
quadruple $(G', X', Y', Z')$, a morphism 
of schemes $X \to X'$, a morphism of schemes
$e : Y' \to Y$ over $X'$, a morphism of group-schemes 
$\phi : G \times X \to G' \times X'$ over $X'$, and a 
$G \times X$ -equivariant morphism of 
schemes $f : Z \to Z'$ over $X'$, the diagram
\[ \begin{CD}
\opn{Hom}_X(\bsym{\Delta}^q_{\K} \times Y, Z)
^{\bsym{\Delta}^q_{0}}
@>{\opn{wsym}_G}>>
\opn{Hom}_X(\bsym{\Delta}^q_{\K} \times Y, Z)
^{\bsym{\Delta}^q_{0}} \\
@V{(e, f)}VV @V{(e, f)}VV \\
\opn{Hom}_{X'}(\bsym{\Delta}^q_{\K} \times Y', Z')
^{\bsym{\Delta}^q_{0}}
@>{\opn{wsym}_{G'}}>>
\opn{Hom}_{X'}(\bsym{\Delta}^q_{\K} \times Y', Z')
^{\bsym{\Delta}^q_{0}} 
\end{CD} \]
is commutative, and likewise for $\opn{w}_G$.
\end{lem}

\begin{proof}
This is due to the functoriality of the exponential map, see Lemma 
\ref{lem1.1}.
\end{proof}

\begin{lem} \label{lem6.X3}
Assume $G$ is abelian. Then for any
$(f_0, \ldots, f_q) \in \opn{Hom}_X(\bsym{\Delta}^q_{\K} \times Y, Z)
^{\bsym{\Delta}^q_{0}}$
the sequence
$\opn{wsym}_{G}(f_0, \ldots, f_q)$ is constant. 
\end{lem}

\begin{proof}
In this case $\opn{exp} : \mfrak{g} \to G$ is an isomorphism of 
algebraic groups, where $\mfrak{g}$ is viewed as an additive 
group. So we may assume that $Z$ is a torsor under 
$\mfrak{g} \times X$. Let
$(f'_0, \ldots, f'_q) := \opn{wsym}_{G}(f_0, \ldots, f_q)$,
and let 
$\gamma_{i,j} : \bsym{\Delta}^q_{\K} \times Y \to \mfrak{g}$
be morphisms such that 
$f_j  = \gamma_{i, j} + f_i$. 
Take $(w, y) \in \bsym{\Delta}^q_{\K} \times Y$. Then 
\[ f'_i(w,y) = 
\Big( \sum_{j=0}^q t_j(w) \cdot \gamma_{i,j}(w,y) \Big)
+ f_i(w,y) \]
for any $i$. Because
$\gamma_{i,j} = - \gamma_{j,i} = \gamma_{0,j} - \gamma_{0,i}$,
$f_i = f_0 + \gamma_{0,i}$ and 
$\sum_{j=0}^q t_j(w) = 1$, 
it follows that
$f'_i(w,y) = f'_0(w,y)$.
\end{proof}

Let's write $\opn{wsym}_{G}^d$ for the $d$-th iteration of the 
operation $\opn{wsym}_{G}$.

\begin{lem} \label{lem6.X4}
For any
$\bsym{f} = (f_0, \ldots, f_q) 
\in \opn{Hom}_X(\bsym{\Delta}^q_{\K} \times Y, Z)
^{\bsym{\Delta}^q_{0}}$
the sequence $\opn{wsym}_{G}^{d(G)}(\bsym{f})$ is constant.
For any $d \geq d(G)$ one has
$\opn{wsym}_{G}^{d}(\bsym{f}) = \opn{wsym}_{G}^{d(G)}(\bsym{f})$.
\end{lem}

\begin{proof}
For any $k$, the orbit of 
$f_0 \in \opn{Hom}_X(\bsym{\Delta}^q_{\K} \times Y, Z)$ 
under the action of the group
$G_k(\bsym{\Delta}^q_{\K} \times Y)$ will be denoted by
$G_k(\bsym{\Delta}^q_{\K} \times Y) \cdot f_0$.
Let 
$(f'_0, \ldots, f'_q) := \opn{wsym}_{G}(f_0, \ldots, f_q)$.
We will prove that if
$f_1, \ldots, f_q \in G_k(\bsym{\Delta}^q_{\K} \times Y) \cdot 
f_0$
then 
$f'_1, \ldots, f'_q
\in G_{k+1}(\bsym{\Delta}^q_{\K} \times Y) \cdot f'_0$.
The assertions of the lemma will then follow. 

Let $\til{Y} = \til{X} :=  \bsym{\Delta}^q_{\K} \times Y$ and
$\til{Z} := \til{X} \times_X Z$. 
So $\til{Z}$ is a torsor under $G \times \til{X}$, and $f_0$
induces a morphism 
$\til{f}_0 \in \opn{Hom}_{\til{X}}(\til{Y}, \til{Z})$.
The morphism
$\tau : G \times \til{X} \to \til{Z}$, 
$(g, \til{x}) \mapsto g \cdot \til{f}_0(\til{x})$,
is an isomorphism of $\til{X}$-schemes. Define 
$\til{W} := \tau(G_k \times \til{X}) \subset \til{Z}$. 
Then $\til{W}$ 
is the ``geometric orbit'' of $\til{f}_0$ under 
$G_k \times \til{X}$; and in particular $\til{W}$ is a torsor under 
$G_k \times \til{X}$. By assumption 
$\til{f}_1, \ldots, \til{f}_q \in 
\opn{Hom}_{\til{X}}(\til{Y}, \til{W})$. 
Define
$(\til{f}'_0, \ldots, \til{f}'_q) :=
\opn{wsym}_{G_k}(\til{f}_0, \ldots, \til{f}_q)$.
By Lemma \ref{lem6.X2} it suffices to prove that 
$\til{f}'_1, \ldots, \til{f}'_q \in
G_{k+1}(\til{Y}) \cdot \til{f}'_0$.

Define $\bar{W} := \til{W} / G_{k+1}$. 
This is a torsor under the group scheme
$(G_k / G_{k+1}) \times \til{X}$. Let
$\bar{f}_0, \ldots, \bar{f}_q \in 
\opn{Hom}_{\til{X}}(\til{Y}, \bar{W})$
be the images of $(\til{f}_0, \ldots, \til{f}_q)$.
Because the group $G_k / G_{k+1}$ is abelian, Lemma \ref{lem6.X3} 
says that 
$\opn{wsym}_{G_k / G_{k+1}}(\bar{f}_0, \ldots, \bar{f}_q)$
is a constant sequence. 
Again using Lemma \ref{lem6.X2}, we see that in fact
$\til{f}'_1, \ldots, \til{f}'_q \in
G_{k+1}(\til{Y}) \cdot \til{f}'_0$.
\end{proof}

Given an $X$-scheme $Y$ the collections
$\bigl\{ \opn{Hom}_{X}(\bsym{\Delta}^q_{\K} 
\times Y, Z) \bigr\}_{q \in \mbb{N}}$ and \linebreak
$\bigl\{ \opn{Hom}_{X}(Y, Z)
^{\bsym{\Delta}^q_0} \bigr\}_{q \in \mbb{N}}$
are simplicial sets. For $q = 0$ there are equalities
\begin{equation} \label{eqn6.2}
\opn{Hom}_{X}(\bsym{\Delta}^0_{\K} \times Y, Z) = 
\opn{Hom}_{X}(Y, Z) 
= \opn{Hom}_{X}(Y, Z)^{\bsym{\Delta}^0_0} . 
\end{equation}

\begin{thm} \label{thm5.1}
Let $G$ be a unipotent algebraic group over 
$\K$, let $X$ be a $\K$-scheme, and let $Z \to X$ be a $G$-torsor 
over $X$. For any $X$-scheme $Y$ and natural number $q$ 
there is a function
\[ \opn{wav}_G : \opn{Hom}_X(Y, Z)^{\bsym{\Delta}^q_{0}} \to
\opn{Hom}_X(\bsym{\Delta}^q_{\K} \times Y, Z) \]
called the {\em weighted average}. The function $\opn{wav}_G$ 
enjoys the following properties.
\begin{enumerate}
\item Symmetric: $\opn{wav}_G$ is equivariant for the  
action of the permutation group of $\{ 0, \ldots, q \}$ on 
$\bsym{\Delta}^q_{0}$ and on $\bsym{\Delta}^q_{\K}$.
\item Simplicial: $\opn{wav}_G$ is a map of simplicial sets 
\[ \bigl\{ \opn{Hom}_{X}(Y, Z)
^{\bsym{\Delta}^q_0} \bigr\}_{q \in \mbb{N}}  \to
\bigl\{ \opn{Hom}_{X}(\bsym{\Delta}^q_{\K} 
\times Y, Z) \bigr\}_{q \in \mbb{N}} . \]
\item Functorial: given another such 
quadruple $(G', X', Y', Z')$, a morphism 
of schemes $X \to X'$, a morphism of $X'$-group-schemes 
$G \times X \to G' \times X'$,
a $G \times X$ -equivariant morphism of 
$X'$-schemes $f : Z \to Z'$ and a morphism of $X'$-schemes
$e : Y' \to Y$, the diagram
\[ \begin{CD}
\opn{Hom}_X(Y, Z)^{\bsym{\Delta}^q_{0}}
@>{\opn{wav}_G}>>
\opn{Hom}_X(\bsym{\Delta}^q_{\K} \times Y, Z) \\
@V{(e, f)}VV @V{(e, f)}VV \\
\opn{Hom}_{X'}(Y', Z')^{\bsym{\Delta}^q_{0}}
@>{\opn{wsym}_{G'}}>>
\opn{Hom}_{X'}(\bsym{\Delta}^q_{\K} \times Y', Z')
\end{CD} \]
is commutative.
\item If $q = 0$ then $\opn{wav}_G$ is the identity map of 
$\opn{Hom}_X(Y, Z)$.
\end{enumerate}
\end{thm}

\begin{proof}
Given a sequence
$\bsym{f} = 
(f_0, \ldots, f_q) \in \opn{Hom}_X(Y, Z)^{\bsym{\Delta}^q_{0}}$
define 
$\opn{wav}_G(\bsym{f}) := f' \in  
\opn{Hom}_X(\bsym{\Delta}^q_{\K} \times Y, Z)$
to be the morphism such that
\[ (\opn{wsym}^{d(G)}_G \circ \opn{w}_G)(f_0, \ldots, f_q) = 
(f', \ldots, f') ; \]
see Lemma \ref{lem6.X4}. Properties (1)-(4) follow from the 
corresponding properties of $\opn{w}_G$ and $\opn{wsym}_G$.
\end{proof}

\begin{proof}[Proof of Corollary \tup{\ref{cor0.1}}]
Take $X = Y := \opn{Spec} \K$ in Theorem \ref{thm5.1}, and 
consider the $G$-torsor $\ul{Z} := G$. 
Choose any base point $z \in Z$;
this defines an isomorphism of left $G(\K)$-sets 
$\ul{Z}(\K) \cong Z$.
The weight sequence $\bsym{w}$ can be considered as a 
$\K$-rational point of $\bsym{\Delta}^q_{\K}$, and we define
\[ \opn{wav}_{G, \bsym{w}}(\bsym{z}) :=
\opn{wav}_{G}(\bsym{z})(\bsym{w}) \in Z . \]
If we were to choose another base point $z' \in Z$ this would 
amount to applying an automorphism of the torsor $\ul{Z}$,
namely right multiplication by some element of $G(\K)$. Due to the 
functoriality of $\opn{wav}_{G}$ the point  
$\opn{wav}_{G, \bsym{w}}(\bsym{z})$ will be unchanged.

The properties of this set-theoretical
averaging process are now immediate consequences of 
the corresponding properties of the geometric average. 
\end{proof}

\begin{rem} \label{rem1.1}
Z. Reichstein observed that our averaging process provides a new 
proof (in characteristic $0$)
of the fact that a unipotent group $G$ is special, namely 
any $G$-torsor $Z$ over $\K$ has a $\K$-rational point. Let us 
explain the idea. 

Let $z_0 \in Z$ be some closed point. Choose a finite Galois 
extension $L$ of $\K$ containing the residue field 
$\bsym{k}(z_0)$. Let $\Gamma$ be the Galois group of $L$ over 
$\K$, which acts on the set $Z(L)$.
Let $z_0, \ldots, z_q \in Z(L)$ 
be the $\Gamma$-conjugates of $z_0$. The group $\Gamma$ acts on 
the sequence 
$\bsym{z} := (z_0, \ldots, z_q)$ by permutations. 
Thus the simultaneous action of $\Gamma$ on 
\[ Z(L)^{\bsym{\Delta}^q_{0}} = 
\opn{Hom}_{\opn{Spec} \K}(\opn{Spec} L, Z)^{\bsym{\Delta}^q_{0}} \]
fixes $\bsym{z}$.

We know that the operator $\mrm{wav}_G$ is symmetric. 
And functoriality says that the action of the Galois group 
on $\opn{Spec} L$ is also respected. 
Since $\bsym{z}$ is fixed by the simultaneous action of
$\Gamma$, so is $\opn{wav}_G(\bsym{z})$. 
Take the uniform weight sequence 
$\bsym{w} := (\frac{1}{q+1}, \ldots, \frac{1}{q+1})$
and define
$z' := \opn{wav}_G(\bsym{z})(\bsym{w}) \in Z(L)$.
Because $\bsym{w}$ is fixed by the permutation group we conclude 
that $z'$ is $\Gamma$-invariant, and hence $z' \in Z(\K)$.
\end{rem}

\begin{rem}
Theorem \ref{thm5.1} has a rather obvious parallel in differential 
geometry. Indeed, a simply connected nilpotent Lie group is the 
same as the group $G(\mbb{R})$ of rational points of a unipotent 
algebraic group $G$ over $\mbb{R}$. 
\end{rem}

\section{Simplicial Sections}
\label{sec2}

In this section we show how the averaging process is used to 
obtain simplicial sections of certain bundles. 

Suppose $H$ and $G$ are affine group schemes over $\K$, and $H$ 
acts on $G$ by automorphisms. Namely there is a morphism of 
schemes $H \times G \to G$ which  for every $\K$-scheme $Y$ 
induces a group homomorphism
$H(Y) \to \opn{Aut}_{\cat{Groups}}(G(Y))$.
Then $H \times G$ has a structure of a group scheme, and we denote 
this group by $H \ltimes G$; it is a geometric semi-direct 
product. 

Recall that an affine group scheme $G$ is called
{\em pro-unipotent} if it is isomorphic to an inverse limit
$\lim_{\leftarrow i} G_i$ of an inverse system
$\{ G_i \}_{i \geq 0}$ of (finite type affine) unipotent groups. 
One may assume that each of the morphisms
$G \to G_{i} \to G_{i-1}$ is surjective. 
Thus $G_i \cong G / N_i$ where $N_i$ is a 
normal closed subgroup of $G$.

We will be concerned with the following geometric situation.

\begin{sce} \label{sce2.1}
Let $H \ltimes G$ be an affine group scheme over $\K$.
Assume $G$ is pro-unipotent, and moreover there exists a 
sequence $\{ N_i \}_{i \geq 0}$ of $H$-invariant 
closed normal subgroups of $G$ such that 
$G \cong \lim_{\leftarrow i} G / N_i$ and each 
$G / N_i$ is unipotent. Let $\pi : Z \to X$ be an 
$H \ltimes G$ -torsor over $X$ 
which is locally trivial for the Zariski 
topology of $X$. Define $\bar{Z} := Z / H$ and let 
$\bar{\pi} : \bar{Z} \to X$ be the projection.
\end{sce}

\begin{thm} \label{thm5.2}
Assume Scenario \tup{\ref{sce2.1}}.
Suppose $U \subset X$ is an open set 
and \linebreak
$\sigma_0, \ldots, \sigma_q : U \to \bar{Z}$ are sections of 
$\bar{\pi}$. Then there exists a morphism 
\[ \sigma : \bsym{\Delta}^q_{\K} \times U \to \bar{Z} \]
such that the diagram
\[ \begin{CD}
\bsym{\Delta}^q_{\K} \times U @>{\sigma}>>
\bar{Z} \\
@V{\mrm{p}_2}VV @V{\bar{\pi}}VV \\
U @>{}>> X 
\end{CD} \]
is commutative. The morphism $\sigma$ depends functorially 
on $U$ and simplicially on the sequence 
$(\sigma_0, \ldots, \sigma_q)$.
If $q = 0$ then $\sigma = \sigma_0$.
\end{thm}

\begin{proof}
We might as well assume that $U = X$. Consider the quotient 
$Z / G$. Since $G$ is normal in $H \ltimes G$ it follows that
$Z / G$ is a torsor under $H \times X$. Let's write 
$\pi_H : Z \to \bar{Z} = Z / H$ and 
$\pi_G : Z \to Z / G$ for the projections. 

Pick an open set $V \subset X$ which trivializes $\pi : Z \to X$. 
Let's write $Z|_V := \pi^{-1}(V)$. 
Because $\pi_H|_V : Z|_V \to \bar{Z}|_V$ is a trivial torsor under
$H \times \bar{Z}|_V$, we can lift the sections 
$\sigma_0, \ldots, \sigma_q$ to sections
$\til{\sigma}_0, \ldots, \til{\sigma}_q : V \to Z$
such that $\pi_H \circ \til{\sigma}_j  = \sigma_j$.
Furthermore, since $\pi_G : Z \to Z / G$ is $H$-equivariant
and $Z / G$ is a torsor under $H \times X$, it follows that we can 
choose $\til{\sigma}_0, \ldots, \til{\sigma}_q$ such that 
$\pi_G \circ \til{\sigma}_j = \tau$ for some 
section $\tau : V \to Z / G$.

Let $F \subset Z|_V$ be the fiber over $\tau$, i.e.\
$F := V \times_{Z / G} Z$ via the morphisms 
$\pi_G : Z \to Z / G$ and $\tau : V \to Z / G$. Then $F$ is a 
torsor under $G \times V$, and
$\til{\sigma}_0, \ldots, \til{\sigma}_q \in 
\opn{Hom}_X(V, F)$.
See diagram below. 
\[ \UseTips  \xymatrix @C=5ex @R=5ex {
& Z 
\ar[dl]_{\pi_H}
\ar[dd]_{\pi}
\ar[dr]^{\pi_G}
& & F 
\ar[ll]_{\supset} \\
\bar{Z} = Z / H 
\ar[dr]^{\bar{\pi}}
& & Z / G
\ar[dl] \\
& X \\
& & & V
\ar[ull]_{\supset}
\ar@{-->}[uul]^{\tau}
\ar@{-->}[uuu]^{\til{\sigma}_i}
\ar@/^5.5ex/[llluu]^{\sigma_i}
} \]

For any $i$ define $F_i := F / N_i$, which is a torsor 
under $(G / N_i) \times V$. Let $\alpha_i : F \to F_i$ be the 
projection, so 
$\alpha_i \circ \til{\sigma}_j \in \opn{Hom}_{X}(V, F_i)$.
By Theorem \ref{thm5.1} we get an average
\begin{equation} \label{eqn2.6}
\rho_i := \opn{wav}_{G / N_i}(\alpha_i \circ \til{\sigma}_0, 
\ldots, \alpha_i \circ \til{\sigma}_q) : 
\bsym{\Delta}^q_{\K} \times V \to F_i .
\end{equation}
The functoriality of $\opn{wav}$ says that the $\rho_i$ form an 
inverse system, and we let
\begin{equation} \label{eqn2.7}
\rho := \lim_{\leftarrow i} \rho_i : 
\bsym{\Delta}^q_{\K} \times V \to F 
\end{equation}
and
\begin{equation} \label{eqn5.5}
\sigma := \pi_H \circ \rho : 
\bsym{\Delta}^q_{\K} \times V \to \bar{Z} . 
\end{equation}

We claim that the morphism $\sigma$ does not depend 
on the choice of the section $\tau : V \to Z / G$. 
Suppose $\tau' : V \to Z / G$ is another such section.
Let $F'$ be the fiber over $\tau'$, and let
$\rho' : \bsym{\Delta}^q_{\K} \times V \to F'$
be the corresponding morphism as in (\ref{eqn2.7}). Now
$\tau' = h \cdot \tau$ for some morphism $h : V \to H$. 
Then $F' = h \cdot F$, and $h : F \to F'$ is a $G \times V$
-equivariant morphism of torsors, 
with respect to the group-scheme automorphism
$\opn{Ad}(h) : G \times V \to G \times V$. 
The new lift of $\sigma_j$ is 
$\til{\sigma}'_j := h \cdot \til{\sigma}_j : V \to F'$.
Define $F_i' := F' / N_i$,
and let $\rho'_i : \bsym{\Delta}^q_{\K} \times V \to F'_i$
be the morphism as in (\ref{eqn2.6}). 
Since $N_i \times V = \opn{Ad}(h)(N_i \times V)$,
we get a group-scheme automorphism
$\opn{Ad}(h) : (G/N_i) \times V \to (G/N_i) \times V$,
and a $(G/N_i) \times V$ -equivariant morphism of 
torsors $h : F_i \to h \cdot F'_i$. 
By functoriality of $\opn{wav}$ (property 3 in Theorem \ref{thm5.1})
it follows that $\rho'_i = h \cdot \rho_i$. Therefore 
$\rho' = h \cdot \rho$, and
$\pi_{H} \circ \rho' = \pi_{H} \circ \rho = \sigma$. 

Property 2 in Theorem \ref{thm5.1}
implies that $\sigma$ depends 
simplicially on $(\sigma_0, \ldots, \sigma_q)$.

Finally take an open covering $X = \bigcup V_j$ such that each 
$V_j$ trivializes $\pi : Z \to X$, and let
$\sigma_j : \bsym{\Delta}^q_{\K} \times V_j \to \bar{Z}|_{V_j}$ 
be the morphism constructed in 
(\ref{eqn5.5}). Since no choices were made we have
$\sigma_j|_{V_j \cap V_k} = \sigma_k|_{V_j \cap V_k}$
for any two indices. Therefore these sections can be glued to a 
morphism 
$\sigma : \bsym{\Delta}^q_{\K} \times X \to \bar{Z}$. 
The functorial and simplicial 
properties of $\sigma$ are clear from its construction.
\end{proof}

Let $X$ be a $\K$-scheme, and let
$X = \bigcup_{i = 0}^m U_{(i)}$ be an open covering, with 
inclusions $g_{(i)} : U_{(i)} \to X$. We denote this covering by 
$\bsym{U}$. For any multi-index
$\bsym{i} = (i_0, \ldots, i_q) \in \bsym{\Delta}^m_q$
we write $U_{\bsym{i}} := \bigcap_{j = 0}^q U_{(i_j)}$,
and we define the scheme
$U_{q} := \coprod_{\bsym{i} \in \bsym{\Delta}^m_q}
U_{\bsym{i}}.$
Given $\alpha \in \bsym{\Delta}_p^q$ and 
$\bsym{i} \in \bsym{\Delta}^m_q$ there is an inclusion of open 
sets
$\alpha_* : U_{\bsym{i}} \to U_{\alpha_*(\bsym{i})}$. 
These patch to a morphism of schemes
$\alpha_* : U_q \to U_p$,
making $\{ U_q \}_{q \in \mbb{N}}$ into a simplicial 
scheme. The inclusions $g_{(i)} : U_{(i)} \to X$ induce inclusions
$g_{\bsym{i}} : U_{\bsym{i}} \to X$
and morphisms 
$g_q : U_q \to X$; and one has the relations
$g_p \circ \alpha_* = g_q$
for any $\alpha \in \bsym{\Delta}_p^q$.

\begin{dfn} \label{dfn3.2}
Let $\pi : Z \to X$ be a morphism of $\K$-schemes.
A {\em simplicial section of 
$\pi$} based on the covering $\bsym{U}$ is a sequence of morphisms
\[ \bsym{\sigma} = \{ \sigma_q : 
\bsym{\Delta}^q_{\K} \times U_q \to Z \}_{q \in \mbb{N}}  \]
satisfying the following conditions.
\begin{enumerate}
\rmitem{i} For any $q$ the diagram
\[ \begin{CD}
\bsym{\Delta}^q_{\K} \times U_q @>{\sigma_q}>>
Z \\
@V{\mrm{p}_2}VV @V{\pi}VV \\
U_q @>{g_q}>> X 
\end{CD} \]
is commutative.
\rmitem{ii} For any $\alpha \in \bsym{\Delta}_p^q$
the diagram
\[ \UseTips \xymatrix @=4.5ex {
& \bsym{\Delta}^{p}_{\K} \times U_p
\ar[dr]^{\sigma_{p}} 
\\
\bsym{\Delta}^{p}_{\K} \times U_q
\ar[ur]^{\bsym{1} \times \alpha^*} 
\ar[dr]_{\alpha_* \times \bsym{1}} 
& & Z \\
& \bsym{\Delta}^{q}_{\K} \times U_q
\ar[ur]_{\sigma_{q}} 
} \]
is commutative.
\end{enumerate}
\end{dfn}

\begin{cor} \label{cor2.1}
Assume Scenario \tup{\ref{sce2.1}}.
Let $\bsym{U} = \{ U_{(i)} \}_{i=0}^m$ be an open covering
of $X$. Suppose that for any $i \in \{ 0, \ldots, m \}$ we are 
given some section 
$\sigma_{(i)} : U_{(i)} \to \bar{Z}$ of 
$\bar{\pi}$. Then there exists a simplicial section
\[ \bsym{\sigma} = \{ \sigma_q : 
\bsym{\Delta}^q_{\K} \times U_q \to \bar{Z} \}_{q \in \mbb{N}} \]
based on $\bsym{U}$, such that 
$\sigma_0|_{U_{(i)}} = \sigma_{(i)}$
for all $i \in \{ 0, \ldots, m \}$. 
\end{cor}

\begin{proof}
For any multi-index
$\bsym{i} = (i_0, \ldots, i_q)$ we have sections
$\sigma_{(i_0)}, \ldots, \sigma_{(i_q)} : U_{\bsym{i}} \to 
\bar{Z}$. 
Let
$\sigma_{\bsym{i}} : \bsym{\Delta}^q_{\K} \times U_{\bsym{i}} \to  
\bar{Z}$
be the morphism provided by Theorem \ref{thm5.2}. For fixed $q$ 
these patch to a morphism
$\sigma_q : \bsym{\Delta}^q_{\K} \times U_q \to \bar{Z}$.
The functorial and simplicial properties in Theorem \ref{thm5.2} 
imply that this is a simplicial section.
\end{proof}

\begin{figure}
\includegraphics{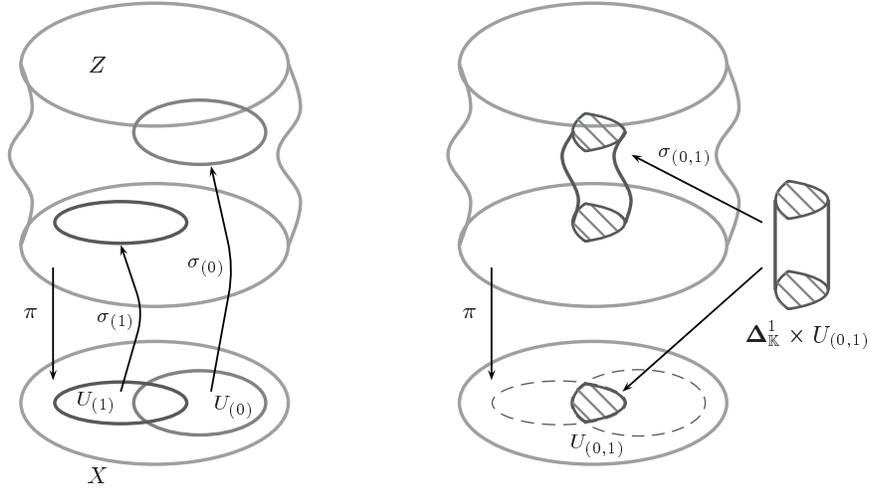}
\caption{Simplicial sections, $q = 1$. 
We start with sections over two open sets $U_{(0)}$ and $U_{(1)}$ 
in the left diagram; and we pass to a simplicial 
section $\sigma_{(0, 1)}$ on the right. As can be seen, 
$\sigma_{(0, 1)}$ interpolates between $\sigma_{(0)}$ and 
$\sigma_{(1)}$.}
\label{fig1}
\end{figure}

This result (with $H$ trivial) is illustrated in Figure \ref{fig1}.

\end{document}